\documentclass[]{amsart}

\usepackage{bbding}
\usepackage[dvips]{graphicx}
\usepackage{amsfonts,amssymb,amsmath,amsthm,amscd}
\allowdisplaybreaks
\usepackage{caption}
\usepackage{subcaption}
\usepackage{mathrsfs}
\usepackage{xcolor}
\usepackage{empheq}
\usepackage{adforn}
\usepackage{cancel}
\usepackage{xr}
\usepackage{mdframed}
\usepackage{tikz-cd}
\usepackage{tikz}
\usepackage{amsfonts}
\usepackage{mathdots}
\usepackage{pinlabel}

\usepackage{enumerate}
\usepackage{lmodern}
\usepackage{mdframed}
\usepackage{stmaryrd}
\usepackage{blindtext}
\usepackage{leftidx}
\usepackage{easytable}
\usepackage{enumitem}
\usepackage{marginnote} 
\usepackage{easyReview}

\usepackage{tgpagella}
\usepackage[T1]{fontenc}

\usepackage[marginparwidth=2.5cm]{geometry}
\definecolor{navie}{RGB}{0, 70, 140}
\usepackage[bookmarks=true,colorlinks=true,urlcolor=red!60!black,linkcolor=red!60!black,citecolor=navie,pdfstartview=FitV]{hyperref}

\makeatletter
\newtheorem*{rep@theorem}{\rep@title}
\newcommand{\newreptheorem}[2]{%
	\newenvironment{rep#1}[1]{%
		\def\rep@title{#2 \ref{##1}}%
		\begin{rep@theorem}}%
		{\end{rep@theorem}}}
\makeatother

\makeatletter
\newtheorem*{rep@cor}{\rep@title}
\newcommand{\newrepcor}[2]{%
	\newenvironment{rep#1}[1]{%
		\def\rep@title{#2 \ref{##1}}%
		\begin{rep@cor}}%
		{\end{rep@cor}}}
\makeatother

\makeatletter
\newtheorem*{rep@prop}{\rep@title}
\newcommand{\newrepprop}[2]{%
	\newenvironment{rep#1}[1]{%
		\def\rep@title{#2 \ref{##1}}%
		\begin{rep@prop}}%
		{\end{rep@prop}}}
\makeatother

\newtheorem{corollary}{Corollary}

\newtheorem{theorem}[corollary]{Theorem}

\newtheorem{proposition}[corollary]{Proposition}

\newrepcor{cor}{Corollary}
\newreptheorem{theorem}{Theorem}
\newrepprop{prop}{Proposition}

\newtheorem*{theorem*}{Theorem}
\newtheorem{lemma}[corollary]{Lemma}

\theoremstyle{definition} 
\newtheorem{definition}[corollary]{Definition}

\theoremstyle{remark} 
\newtheorem{remark}[corollary]{Remark}
 
\newtheorem*{remark*}{Remark}

\newcommand{\R}{\mathbb{R}}

\newcommand{\CP}{\mathbb{CP}}
\renewcommand{\H}{\mathbb{H}}

\newcommand{\dif}{\mathsf{d}}

\newcommand{\pa}{\partial}

\begin{document}

\title{On circle patterns and spherical conical metrics}

\author{Xin Nie}
\address{
	Shing-Tung Yau Center of Southeast University, Nanjing 210096, China}
\email{nie.hsin@gmail.com}

\maketitle

\begin{abstract}
The Koebe-Andreev-Thurston circle packing theorem, as well as its generalization to circle patterns due to Bobenko and Springborn, holds for Euclidean and hyperbolic metrics possibly with conical singularities, but fails for spherical metrics because of the nonuniqueness coming from M\"obius transformations. In this paper, we show that a unique existence result for circle pattern with spherical conical metric holds if one prescribes the total geodesic curvature of each circle instead of the cone angles.  
\end{abstract}


\section{Introduction}\label{sec_1}
\begin{figure}[ht]
	\centering
\includegraphics[width=7cm]{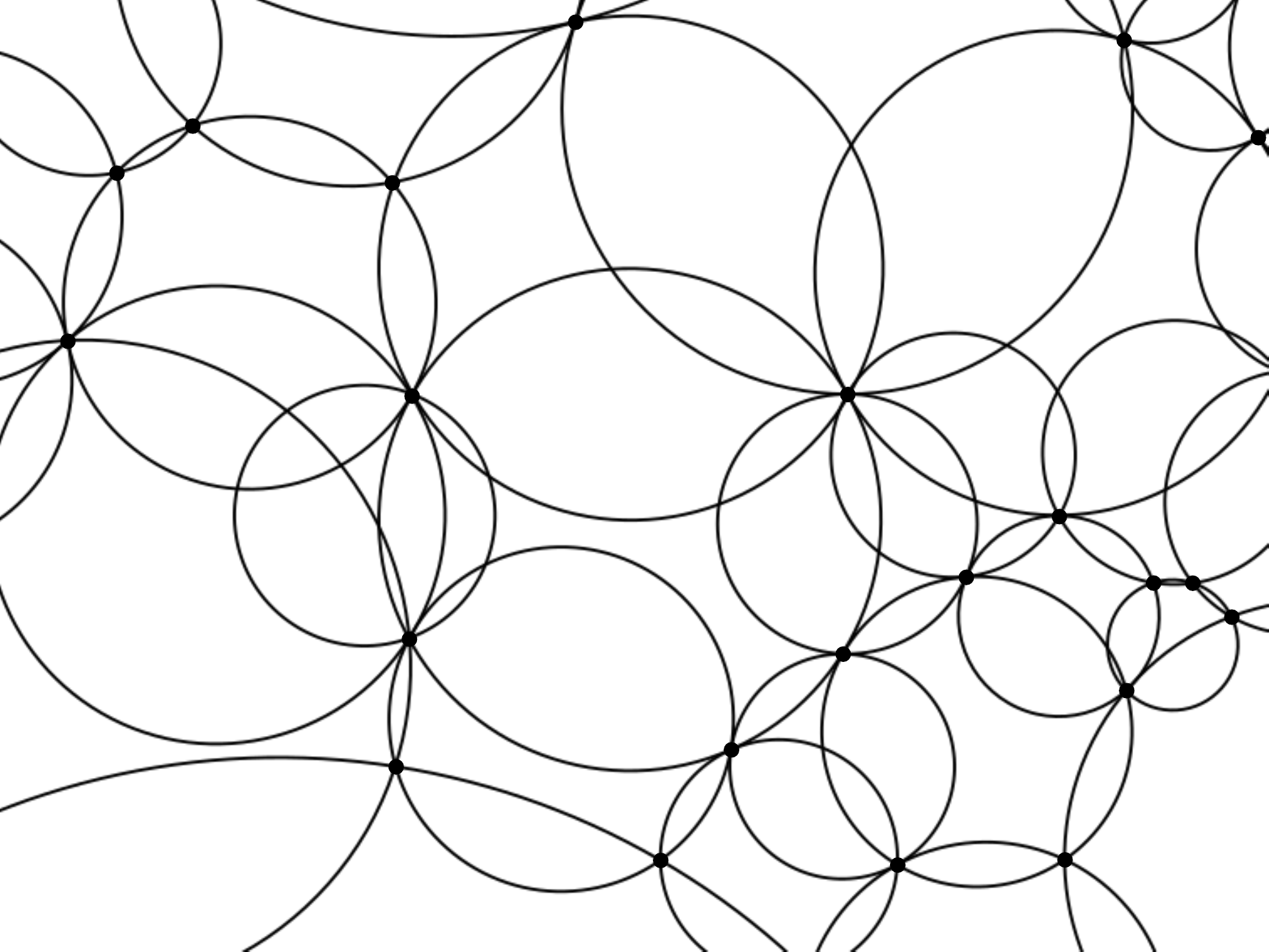}
\hspace{1cm}
\includegraphics[width=7cm]{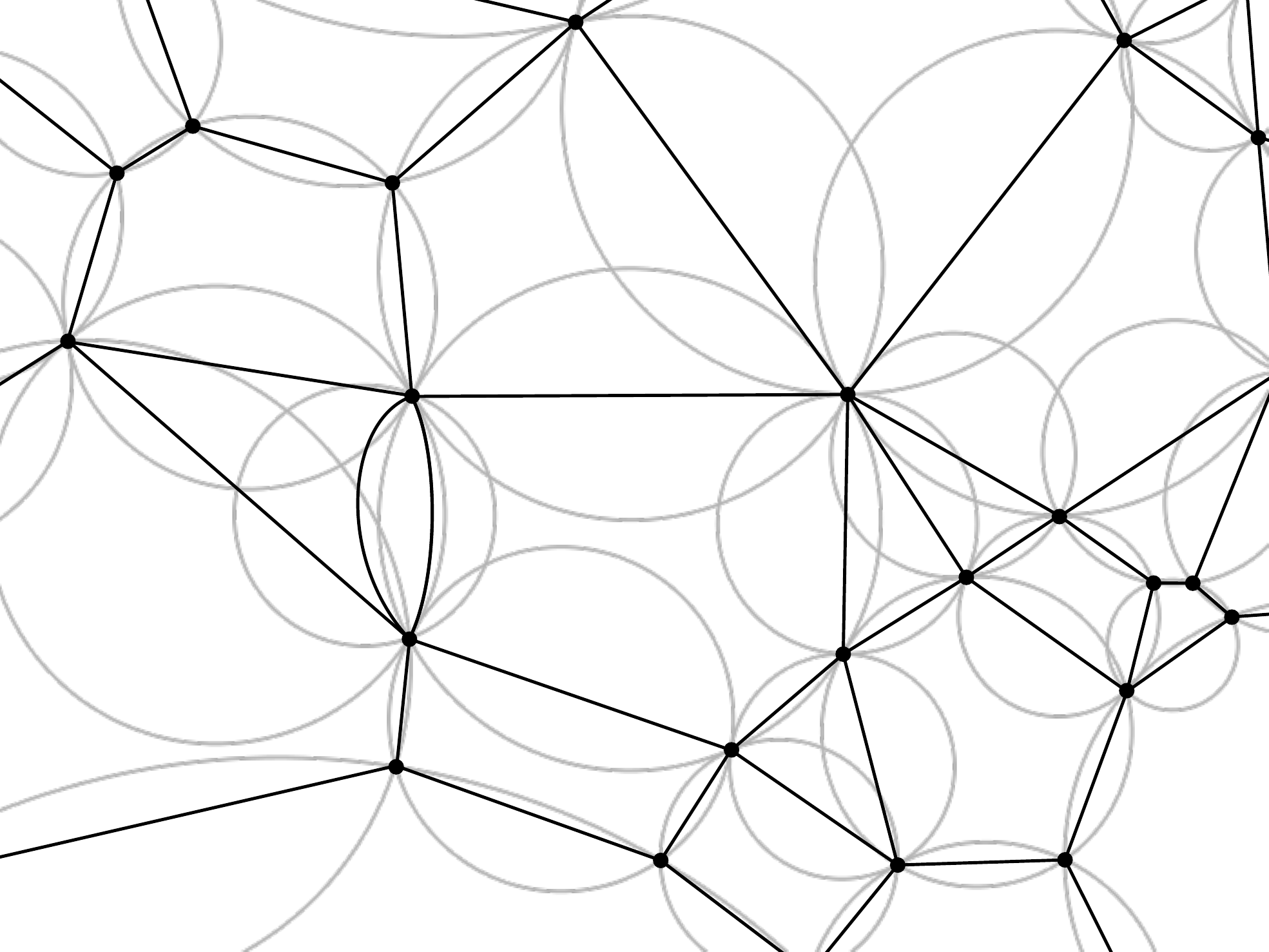}
\caption{A local picture of a circle pattern.}
\label{figure_pattern}
\end{figure}

On a closed surface $\Sigma$ equipped with a spherical, Euclidean or hyperbolic metric, a (Delaunay) circle pattern is a finite set $\mathcal{D}$ of open round disks that looks typically as in Figure \ref{figure_pattern}.  
One associates to $\mathcal{D}$ a weighted graph $\mathsf{G}_\mathcal{D}$ on $\Sigma$ as shown on the right, whose vertices $\mathsf{V}_\mathcal{D}$, edges $\mathsf{E}_\mathcal{D}$ and faces $\mathsf{F}_\mathcal{D}$ are:
\begin{itemize}
	\item 
$\mathsf{V}_\mathcal{D}=\Sigma\setminus\bigcup_{D\in\mathcal{D}}D$;
\item each $f\in \mathsf{F}_\mathcal{D}$ corresponds to a disk $D_f\in\mathcal{D}$;
\item each $e\in \mathsf{E}_\mathcal{D}$ corresponds to a bigon where certain two disks meet, and is weighted by the interior angle $\theta_e\in(0,\pi)$ of that bigon.
\end{itemize}

In this paper, we assume that the above bigons are acute or right angled, namely $\theta_e\in(0,\frac{\pi}{2}]$, but allow $\Sigma$ to have conical singularities at the vertices $\mathsf{V}_\mathcal{D}$ or the centers of the disks. See \S \ref{sec_2} for details. In this setting, the bigons are disjoint from each other (unlike in Figure \ref{figure_pattern}) and do not contain any singularity. The cone angle $\alpha_v$ of $\Sigma$ at a vertex $v\in \mathsf{V}_\mathcal{D}$ is determined by the weights $(\theta_e)$ via the relation
$$
\alpha_v=\sum_{\scalebox{0.65}{\text{$e$: edge from $v$}}}(\pi-\theta_e)~;
$$  
while the cone angle $\alpha_f$ at the center of $D_f$ also satisfies  certain condition in terms of $(\theta_e)$, which in the hyperbolic case reads 
\begin{equation}\label{eqn_angleconstraint}
\sum_{f\in \mathsf{F}'}\alpha_f<\sum_{
\scalebox{0.65}{\parbox{2.9cm}{$e$: edge incident\\ with any face in $\mathsf{F}'$}}
}\hspace{-0.5cm}2\theta_e\quad \text{for any subset of faces $\mathsf{F}'\subset \mathsf{F}$}.
\end{equation}

Improving on Colin de Verdi\`ere's proof \cite{ColindeVerdiere} of  Koebe-Andreev-Thurston Circle Packing Theorem, Bobenko and Springborn \cite{Bobenko-Springborn} showed that there exists a unique hyperbolic $(\Sigma,\mathcal{D})$ realizing a prescribed weighted graph and with prescribed cone angles, as along as the prescription fulfills the above constraints:
\begin{theorem*}[\cite{Bobenko-Springborn}]\label{thm_main}
Let $S$ be a closed topological surface and $\mathsf{G}$ be a graph on $S$ with the sets of edges and faces denoted by $\mathsf{E}$ and $\mathsf{F}$, such that each $e\in \mathsf{E}$ is weighted by some $\theta_e\in(0,\frac{\pi}{2}]$ and every $f\in\mathsf{F}$ is simply connected.
Then $(\alpha_f)_{f\in\mathsf{F}}\in \R_+^\mathsf{F}$ satisfies \eqref{eqn_angleconstraint} if and only if there exists a hyperbolic metric $\sigma$ on $S$ with conical singularities, along with a circle pattern $\mathcal{D}$ on $(S,\sigma)$, such that $\mathsf{G}_\mathcal{D}$ coincides with $\mathsf{G}$ (up to isotopy) as weighted graphs and $\sigma$ has cone angle $\alpha_f$ at the center of the disk $D_f\in\mathcal{D}$ for all $f\in \mathsf{F}$. Moreover, $(\sigma,\mathcal{D})$ is unique up to isotopy if it exists.  
\end{theorem*}
A similar result also holds for Euclidean metrics. When $\alpha_f=2\pi$ for all $f$, we get nonsingular metrics and this implies the Koebe-Andreev-Thurston theorem. For general $(\alpha_f)$, the result can be viewed as a discrete analogue of the McOwen-Troyanov uniformization theorem \cite{McOwen, Troyanov}. 

In the case of spherical metrics, just as the McOwen-Troyanov theorem, such a cone angle prescription result fails because of the nonuniqueness coming from M\"obius transformations, which preserve the angle but distort the metric. The purpose of this paper is to provide the following unique existence result for spherical metrics, whose statement is almost identical to the above theorem except that we now prescribe a more metric related quantity -- the total geodesic curvature $T(\pa D_f)$ of each circle $\pa D_f$ -- instead of the angle $\alpha_f$. Note that we have
$T(\pa D_f)=\alpha_f-\mathrm{Area}(D_f)$ by Gauss-Bonnet.
\begin{theorem}\label{thm_main}
Let $S$ and $\mathsf{G}$ be as in the above theorem. Then $(T_f)_{f\in\mathsf{F}}\in \R_+^\mathsf{F}$ satisfies \eqref{eqn_angleconstraint} (with $\alpha_f$ replaced by $T_f$) if and only if there exists a spherical metric $\sigma$ on $S$ with conical singularities, along with a circle pattern $\mathcal{D}$ on $(S,\sigma)$, such that $\mathsf{G}_\mathcal{D}$ coincides with $\mathsf{G}$ (up to isotopy) as weighted graphs and $\pa D_f$ has total geodesic curvature $T(\pa D_f)=T_f$ for all $f\in \mathsf{F}$. Moreover, $(\sigma,\mathcal{D})$ is unique up to isotopy if it exists.
\end{theorem}
The reason why $(T_f)$ and $(\alpha_f)$ in the two theorems satisfy exactly the same condition is that for a bigon of angle $\theta$ in $\mathbb{S}^2$ or $\H^2$ formed by two circle arcs (see Figure \ref{figure_triangle}), the constraint on the total geodesic curvatures $T_1,T_2$ of the arcs in the spherical case is the same as the constraint on the angles $\alpha_1,\alpha_2$ of the arcs in the hyperbolic case, namely $T_1,T_2>0$, $T_1+T_2<2\theta$ (see Lemma \ref{lemma_T1T2} and Remark \ref{remark_hyperbolicbigon}).

The rest of the paper is devoted to a proof of Theorem \ref{thm_main} through the standard variational method \cite{Bobenko-Springborn,ColindeVerdiere}. The original method for the hyperbolic cone angle prescription problem is based on a convex function in the variables $u_f=-\log\tanh\frac{r_f}{2}$ ($f\in\mathsf{F}$, $r_f$ is the radius of $D_f$) whose Legendre dual variables are the cone angles $\alpha_f$. We use below the counterpart of this function in the spherical setting and show that it is convex in the variables $K_f=\log\cot r_f$, whose Legendre dual variables are the total geodesic curvatures $T(\pa D_f)$. 

\section{Circle patterns on spherical surfaces}\label{sec_2}
Although circle patterns have been extensively studied, it is somewhat technical to give a precise definition of them on surfaces, and in particular to make the relevant notion of ``self-intersection part of a circle'' rigorous. 
Schlenker and Yarmola \cite{Schlenker-Yarmola} gave a definition for closed $\CP^1$-surfaces of genus $\geq1$ by lifting to the universal cover. In this section, we define a restricted version of circle patterns on closed spherical surfaces of any genus with conical singularities. Everything works for Euclidean or hyperbolic surfaces as well.

\begin{figure}[ht]
	\centering
	\includegraphics[width=6cm]{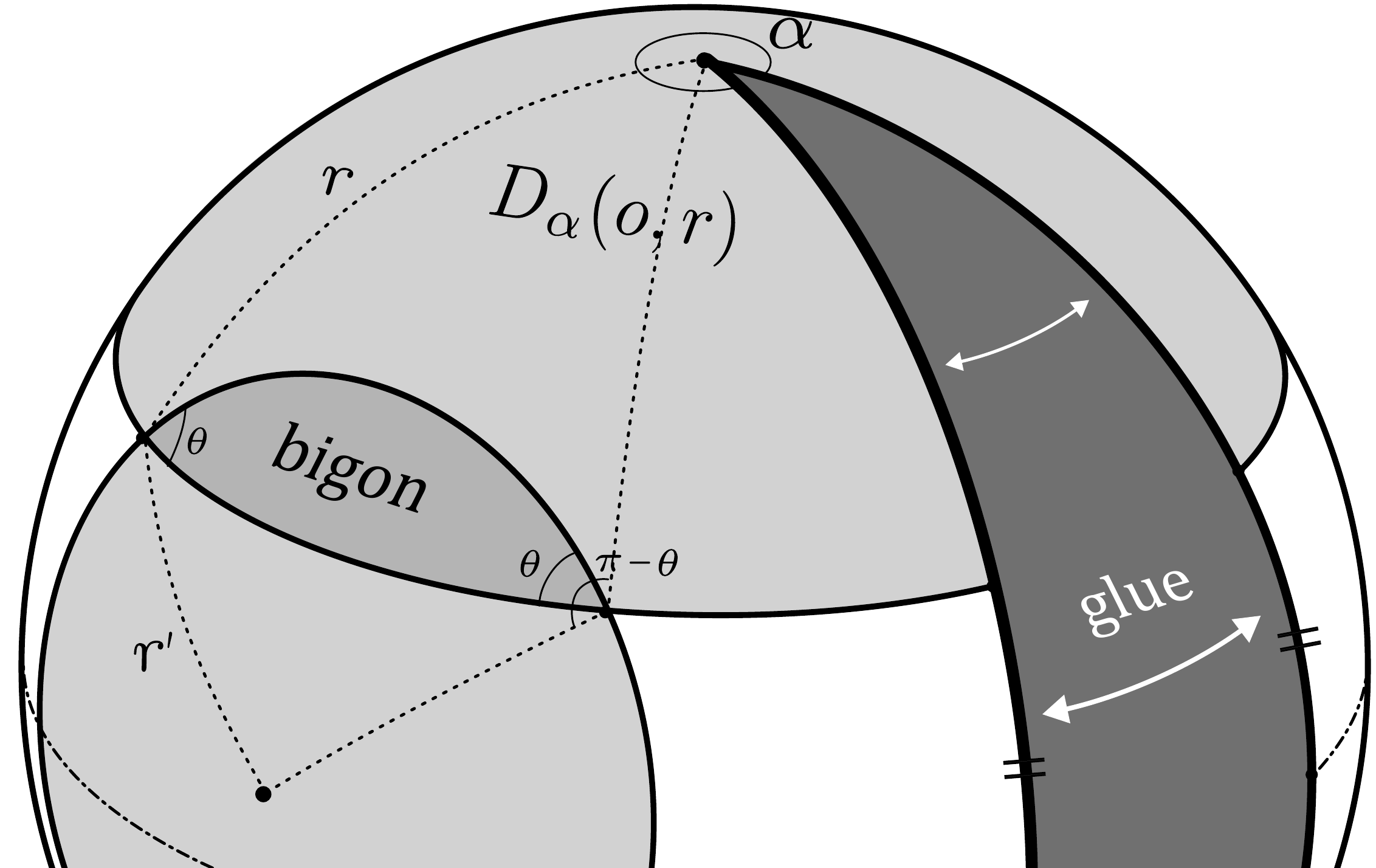}
	\caption{A disk $D_\alpha(o,r)\subset\mathbb{S}^2_\alpha$ and a bigon of angle $\theta$.}
	\label{figure_sphere}
\end{figure}

\begin{figure}[ht]
	\centering
	\includegraphics[width=15cm]{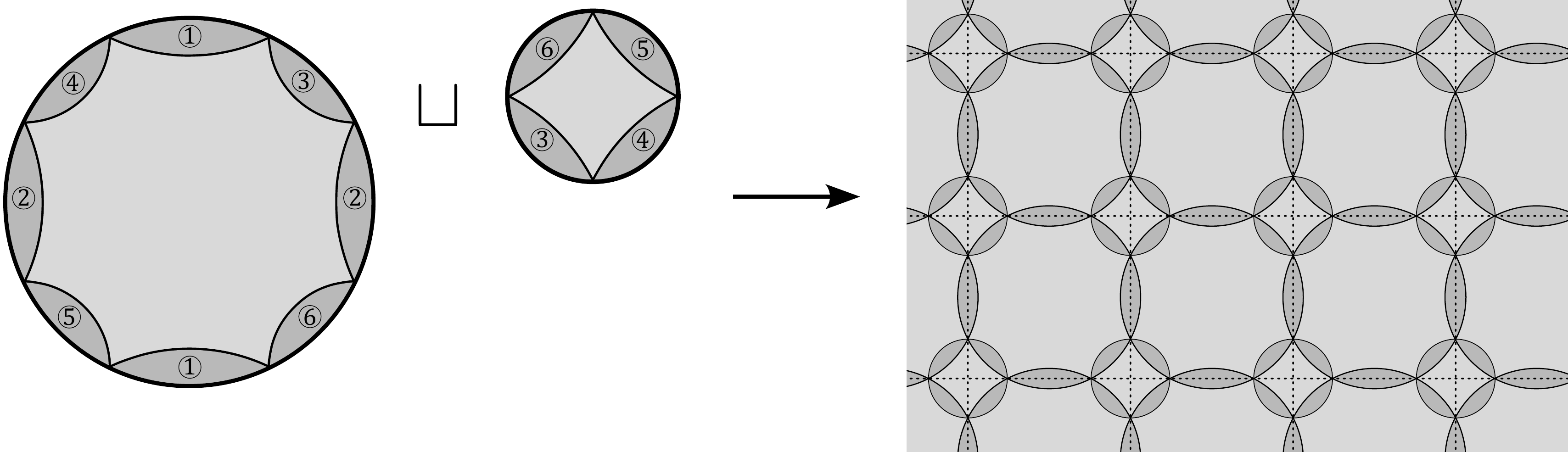}
	\caption{An example of Prop.\ \ref{prop_localiso} with a torus as $\Sigma$ (on the right is the universal cover $\widetilde{\Sigma}$). The pairing of bigons in $\iota^{-1}(B)$ is indicated by labels.}
	\label{figure_pairing}
\end{figure}

A \emph{surface with spherical conical metric}, or a \emph{spherical surface} for short, is a surface locally modeled on the metric spaces $\mathbb{S}^2_\alpha$, $\alpha>0$ shown in Figure \ref{figure_sphere}.
 Each $\mathbb{S}^2_\alpha$ is obtained by gluing together spherical lunes of total angle $\alpha$. Let $D_{\alpha}(o,r)$ denote the open metric disk of radius $r$ in $\mathbb{S}^2_\alpha$ centered at the north pole $o\in\mathbb{S}^2_\alpha$. We will only consider such disks with $r<\frac{\pi}{2}$, or equivalently, the strictly convex ones.
By a \emph{bigon} in a spherical surface, we mean an open set isometric to the intersection of two disks of radii less than $\frac{\pi}{2}$ in $\mathbb{S}^2$ not containing each other. In particular, a bigon does not contain any cone point. The \emph{angle} of a bigon refers to the interior angle $\theta$ form by the two sides (see Figure \ref{figure_sphere}).  
\begin{proposition}\label{prop_localiso}
Take a number of disks $D_i=D_{\alpha_i}(o,r_i)$, $i=1,\cdots,n$ as above, with $\alpha_i\in\R_+$, $r_i\in(0,\frac{\pi}{2})$, and let $\iota$ be a local isometry from the disjoint union $D_1\sqcup\cdots\sqcup D_n$ to a closed spherical surface $\Sigma$ such that 
\begin{itemize}
	\item the image of $\iota$ only misses finitely many points of $\Sigma$;
	\item $\iota$ is at most $2$-to-$1$;
	\item the $\iota(D_i)$'s do not contain each other.
\end{itemize}
Put $B:=\{p\in\Sigma\mid\text{$\iota^{-1}(p)$ has two points}\}$. Then $\iota^{-1}(B)$ 
is a disjoint union of bigons as shown in Figure \ref{figure_pairing}, which fills up the boundary of the $D_i$'s, whereas the equivalence relation ``$\sim$'' on $\iota^{-1}(B)$ defined by $x\sim y$ $\Leftrightarrow$ $\iota(x)=\iota(y)$ is an isometric pairing of these bigons.
\end{proposition}
We omit the proof, which consists of lengthy but elementary arguments using the strict convexity of the disks. We take such $\iota$ as the definition for the restricted version of circle patterns that we study in this paper:
\begin{definition}\label{def_pattern}
A local isometry $\iota:D_1\sqcup\cdots\sqcup D_n\to\Sigma$ as in Proposition \ref{prop_localiso} is called a \emph{circle pattern} on $\Sigma$, and each connected component of the set $B\subset\Sigma$ is called a bigon of the circle pattern. Abusing the notation, we also understand $D_i$ as an immersed disk in $\Sigma$,
and also refer to the set of immersed disks $\mathcal{D}=\{D_1,\cdots,D_n\}$ as a circle pattern (so that each bigon is either an intersection component of two such disk $D_i,D_j$ or a ``self-intersection component'' of a single $D_i$).
Let $\mathsf{G}_\mathcal{D}$ denote the weighted graph on $\Sigma$ such that
\begin{itemize}
	\item the vertices are the finitely many points  $\mathsf{V}_\mathcal{D}:=\Sigma\setminus D_1\cup\cdots\cup D_n$ missed by $\iota$;
	\item the edges $\mathsf{E}_\mathcal{D}$ are in $1$-to-$1$ correspondence with the bigons: given a bigon $b$, the corresponding $e\in\mathsf{E}_\mathcal{D}$ is a curve in $b$ joining the two vertices of $b$;
	\item each $e\in\mathsf{E}_\mathcal{D}$ is weighted by the angle $\theta_e\in(0,\pi)$ of the corresponding bigon.
\end{itemize}
\end{definition}
It is clear from the construction that 
the faces of $\mathsf{G}_\mathcal{D}$ correspond to the disks in $\mathcal{D}$ and are simply connected, hence $\mathsf{G}_\mathcal{D}$ gives a cellular decomposition of $\Sigma$. Conversely, any closed cellular surface can be realized this way by some spherical conical metric and circle pattern. In fact, we show in Proposition \ref{prop_reconstruction} below that there is a unique realization with prescribed angles $(\theta_e)$ in $(0,\frac{\pi}{2}]$ and prescribed radii of disks. 
\begin{remark}
Definition \ref{def_pattern} is more restrictive than the usual notion of Delaunay circle patterns in that the latter allows the bigons to overlap as in Figure \ref{figure_pattern}.
See \cite{Schlenker-Yarmola} for a definition on nonsingular $\CP^1$-surfaces.
Meanwhile, allowing cone points gives our definition more flexibility than the nonsingular setting: since there is no restriction on our $\mathsf{G}_\mathcal{D}$ other than simply connectedness of faces, $\mathsf{G}_\mathcal{D}$ can have a vertex where only one edge is issued, or a face bounded by only one edge, which are both impossible without cone points. 
\end{remark}
\begin{remark}
Given a circle pattern $\mathcal{D}$ on $\Sigma$, the cone points of $\Sigma$ can only be at $\mathsf{V}_\mathcal{D}$ or the centers of the disks. A basic fact, which we do not use explicitly in this paper, is that the cone angle $\alpha_v$ at
$v\in \mathsf{V}_\mathcal{D}$ equals the sum of $\theta'_e:=\pi-\theta_e$ over all edges $e$ issuing from $v$, where $\theta'_e$ appears twice if $e$ is a loop (cf.\ \cite[Lemma 1.11]{Schlenker-Yarmola}).
\end{remark}
From now on, we fix a closed topological surface $S$ and a graph $\mathsf{G}$ as in Theorem \ref{thm_main}, where the faces $\mathsf{F}$ are simply connected and the edges $\mathsf{E}$ are weighted by some $(\theta_e)_{e\in\mathsf{E}}\in (0,\frac{\pi}{2}]^\mathsf{E}$. Consider the moduli space $\mathcal{M}$ of all spherical conical metrics with circle patterns which realize this weighted graph, namely 
$$
\mathcal{M}:=\left\{(\sigma,\mathcal{D})\ \Bigg|\ \parbox{6.1cm}{$\sigma$ is a spherical conical metric on $S$;\\
$\mathcal{D}$ is a circle patter on $(S,\sigma)$ such that $\mathsf{G}_\mathcal{D}$ is isotopic to $\mathsf{G}$ and their weights match
}\right\}/\mathrm{Homeo}^0(S).
$$
By abuse of notation, we sometimes understand a pair $(\sigma,\mathcal{D})$ as the point of $\mathcal{M}$ represented by it. Given $(\sigma,\mathcal{D})$, let $D_f\in\mathcal{D}$ denote the immersed disk in $(S,\sigma)$ corresponding to a face $f\in \mathsf{F}$ and let $r(D_f)\in(0,\frac{\pi}{2})$ be the radius of this disk. By the following proposition, the radius system $(r(D_f))\in(0,\frac{\pi}{2})^\mathsf{F}$
determines $(\sigma,\mathcal{D})$ as a point of $\mathcal{M}$, and hence gives a parametrization $\mathcal{M}\cong(0,\frac{\pi}{2})^\mathsf{F}$:
\begin{proposition}\label{prop_reconstruction}
The map 
$\mathcal{M}\to (0,\tfrac{\pi}{2})^\mathsf{F}$ given by $(\sigma,\mathcal{D})\mapsto (r(D_f))_{f\in \mathsf{F}}$
is bijective.
\end{proposition}
\begin{proof}
The only nontrivial part is the surjectivity. To show this, we construct $(\sigma,\mathcal{D})$ with prescribed radii $(r_f)\in(0,\tfrac{\pi}{2})^\mathsf{F}$ as follows. First, subdivide every face $f\in\mathsf{F}$ by picking a point $p_f\in f$ and adding an edge connecting each vertex on $\pa f$ with $p_f$. This yields a triangulation of $S$. 
The set of triangles, denoted by $\hat{\mathsf{E}}$, identifies with the set of \emph{oriented} edges of $\mathsf{G}$. For any (un-oriented) edge $e\in\mathsf{E}$, letting $e^\pm\in\hat{\mathsf{E}}$ denote the triangles on the two sides of $e$, we endow the quadrilateral $e^+\cup e^-$ with the structure of the specific spherical geodesic quadrilateral shown in Figure \ref{figure_reconstruction},
\begin{figure}[ht]
	\centering
	\includegraphics[width=13.5cm]{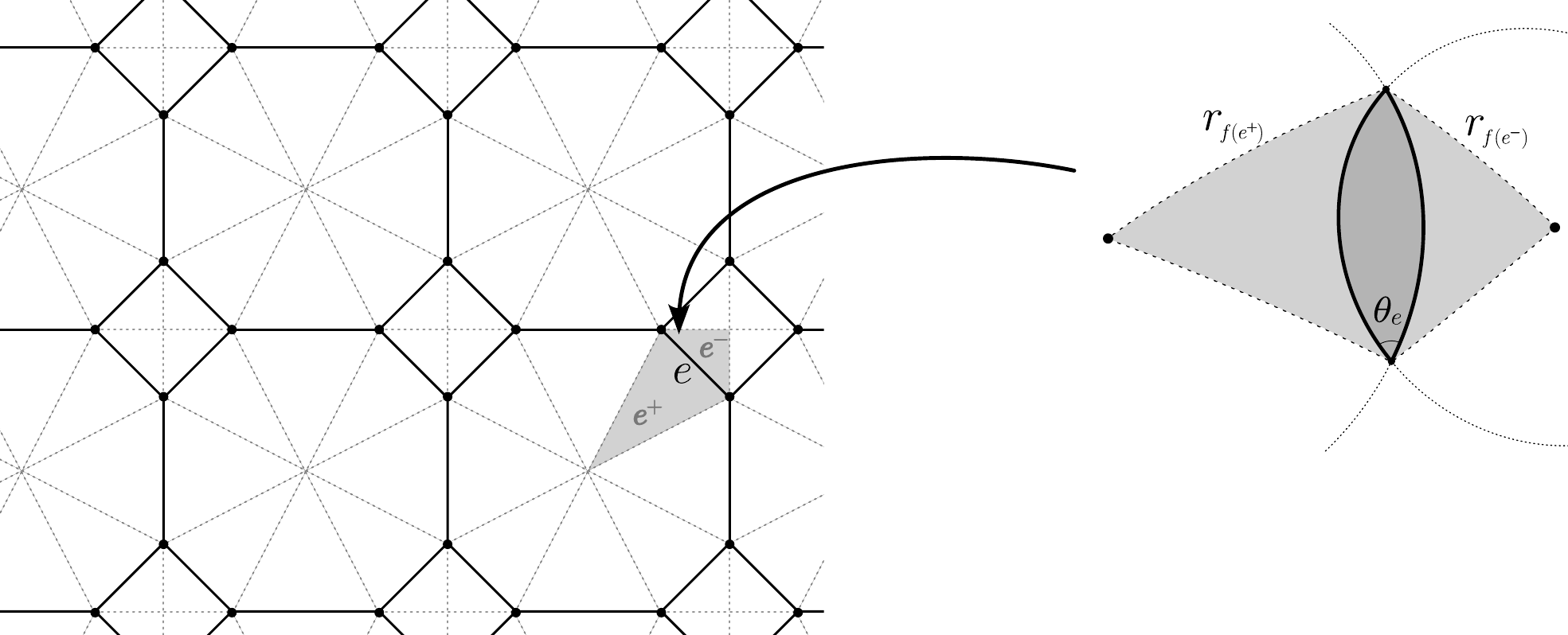}
	\caption{Proof of Prop.\ \ref{prop_reconstruction} (compare Figure \ref{figure_pairing}).}
	\label{figure_reconstruction}
\end{figure} 
 which is the union of two disk sectors with radii $r_{f(e^+)}$ and $r_{f(e^-)}$ respectively, where $f(t)\in \mathsf{F}$ denotes the face containing the triangle $t\in\hat{\mathsf{E}}$. These quadrilaterals fit together and yields a pair  $(\sigma,\mathcal{D})$ with the prescribed radii.
\end{proof}

Since a disk $D\subset\mathbb{S}^2$ of radius $r$ has circumference $2\pi\sin r$ and area $4\pi\sin^2(r/2)=2\pi(1-\cos r)$, by using Gauss-Bonnet, we infer that $\pa D$ has constant geodesic curvature $k(\pa D)=\cot r$. The same holds for any disk $D$ in a circle pattern, regardless of the cone angle at the center.
Instead of the above radius parametrization $\mathcal{M}\cong (0,\tfrac{\pi}{2})^\mathsf{F}$, what we will use below is the reparametrization by the logarithm of geodesic curvature:
$$
\mathcal{M}\overset\sim\to \R^\mathsf{F},\quad (\sigma,\mathcal{D})\mapsto (K_f)_{f\in \mathsf{F}}\ \ \text{with } K_f:=\log k(\pa D_f)=\log \cot r(D_f).
$$
Meanwhile, Theorem \ref{thm_main} claims the unique existence of $(\sigma,\mathcal{D})\in\mathcal{M}$ such that the \emph{total} geodesic curvature $T(\pa D_f)$ of $\pa D_f$, which is by definition $k(\pa D_f)$ times the length $\ell(\pa D_f)$, has the prescribed value. 
We will find $(\sigma,\mathcal{D})$ as the critical point of a function on $\mathcal{M}$ which is strictly convex in the above coordinates $K=(K_f)$. This involves a convex function on the space of spherical bigons discussed in the next section.

\section{A convex function on the space of bigons}\label{sec_3}
Fix $\theta\in(0,\pi)$ throughout this section and let $\mathcal{B}_\theta$ denote the moduli space of spherical bigons of angle $\theta$ with the two sides labeled by 1 and 2. Consider the following geometric quantities of such a bigon:
\begin{itemize}
	\item[$r_i$] \hspace{0.2cm} radius of the disk $D_i\subset\mathbb{S}^2$ whose boundary contains the side $i$ (here and below, $i=1,2$)
	\item[$k_i$] \hspace{0.2cm}geodesic curvature of the side $i$ ($=\cot r_i$)
	\item[$K_i$] \hspace{0.2cm}$:=\log k_i$
	\item[$\alpha_i$] \hspace{0.2cm} angle of the sector of $D_i$ spanned by the side $i$ 
	\item[$\ell_i$] \hspace{0.2cm}length of the side $i$ ($=\alpha_i\sin r_i$)
	\item[$T_i$] \hspace{0.2cm}total geodesic curvature of the side $i$ ($=\ell_ik_i$)
	\item[$A$] \hspace{0.2cm}area ($=2\theta-T_1-T_2$ by Gauss-Bonnet)
\end{itemize}

\begin{figure}[ht]
	\centering
	\includegraphics[width=8cm]{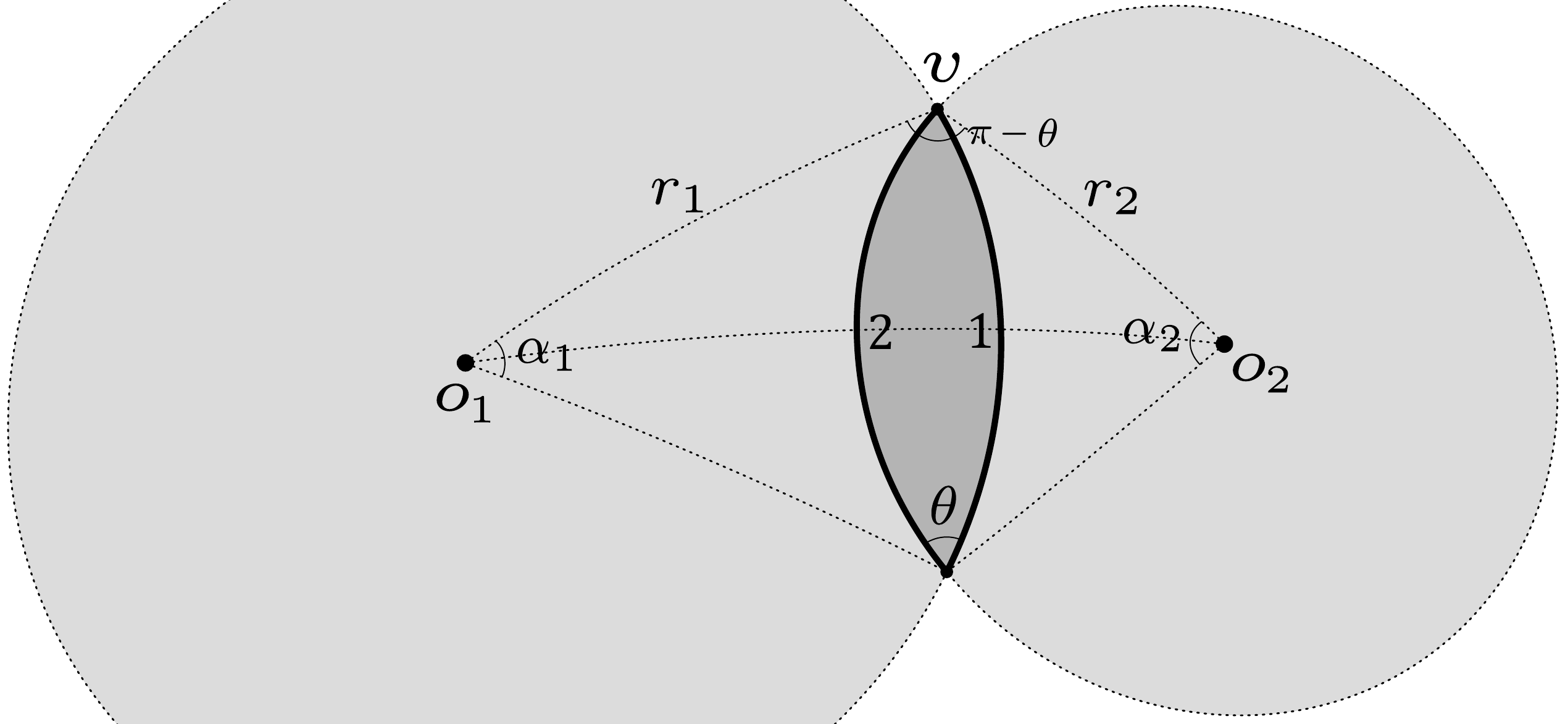}
	\caption{A spherical bigon.}
	\label{figure_triangle}
\end{figure} 

\begin{lemma}\label{lemma_radii}
For any $r_1,r_2\in(0,\frac{\pi}{2})$, there exists, up to isometry, a unique bigon of given angle $\theta$ such that the disk $D_i\subset\mathbb{S}^2$ whose boundary contains the side $i$ has radius $r_i$. 
\end{lemma}
\begin{proof}
An equivalent statement is: for any $0<r_2\leq r_1<\frac{\pi}{2}$, there exists a unique $r_3\in(r_1-r_2,r_1+r_2)$ with the property that if $D_1,D_2\subset\mathbb{S}^2$ are disks with radii $r_1,r_2$ and their centers have distance $r_3$, then the bigon $D_1\cap D_2$ has angle $\theta$. But this follows immediately from the relation 
$$
\cos r_3=\cos r_1\cos r_2+\sin r_1\sin r_2 \cos(\pi-\theta),
$$
which is the cosine law (see \cite{wiki:Spherical_trigonometry}) applied to the triangle $o_1vo_2$ in Figure \ref{figure_triangle}. 
\end{proof}

We now view the above quantities as functions on $\mathcal{B}_\theta$. Then Lemma \ref{lemma_radii} means that
$(r_1,r_2)$ is a parametrization of $\mathcal{B}_\theta$ by $(0,\frac{\pi}{2})^2$. Therefore, $(k_1,k_2)$ and $(K_1,K_2)$ give reparametrizations by $\R_+^2$ and $\R^2$, respectively. 
\begin{lemma}\label{lemma_closed}
The	$1$-form $\omega_\theta:=\ell_1\dif k_1+\ell_2\dif k_2=T_1\dif K_1+T_2\dif K_2$
 on $\mathcal{B}_\theta$ is closed. 
\end{lemma}
\begin{proof}
Viewing $(k_1,k_2)$ as coordinates of $\mathcal{B}_\theta$, we only need to show $\frac{\pa\ell_1}{\pa k_2}=\frac{\pa\ell_2}{\pa k_1}$. To this end, look again at the triangle $o_1vo_2$ in Figure \ref{figure_triangle}, which has angles $\alpha_1':=\frac{\alpha_1}{2},\alpha_2':=\frac{\alpha_2}{2},\theta':=\pi-\theta$ and side lengths $r_1,r_2,r_3$. The \emph{cotangent $4$-part formula} (see \cite{wiki:Spherical_trigonometry}) gives
$\cos r_1\cos\theta'=\cot r_2\sin r_1-\cot\alpha'_1\sin\theta'$, hence
$$
\cot\alpha'_1=\frac{1}{\sin\theta'}\left(\cot r_2\sin r_1-\cos r_1\cos\theta'\right).
$$
Since $\theta'=\text{constant}$, $\ell_1=2\alpha'_1\sin r_1$ and $\cot r_i=k_i$ (in particular, $\frac{\pa r_1}{\pa k_2}=\frac{\pa r_2}{\pa k_1}=0$), by taking the partial derivative in $k_2$ of the above equality, we obtain
$$
\frac{\pa \ell_1}{\pa k_2}=-\frac{2\sin^2 r_1\sin^2\alpha_1'}{\sin\theta'}.
$$
It follows by symmetry that
$$
\frac{\pa \ell_2}{\pa k_1}=-\frac{2\sin^2 r_2\sin^2\alpha_2'}{\sin\theta'}.
$$
By virtue of the sine law $\frac{\sin \alpha_1'}{\sin r_2}=\frac{\sin\alpha_2'}{\sin r_1}$, we concluded that $\frac{\pa\ell_1}{\pa k_2}=\frac{\pa\ell_2}{\pa k_1}$, as required.
\end{proof}

\begin{figure}[ht]
	\centering
	\includegraphics[width=8cm]{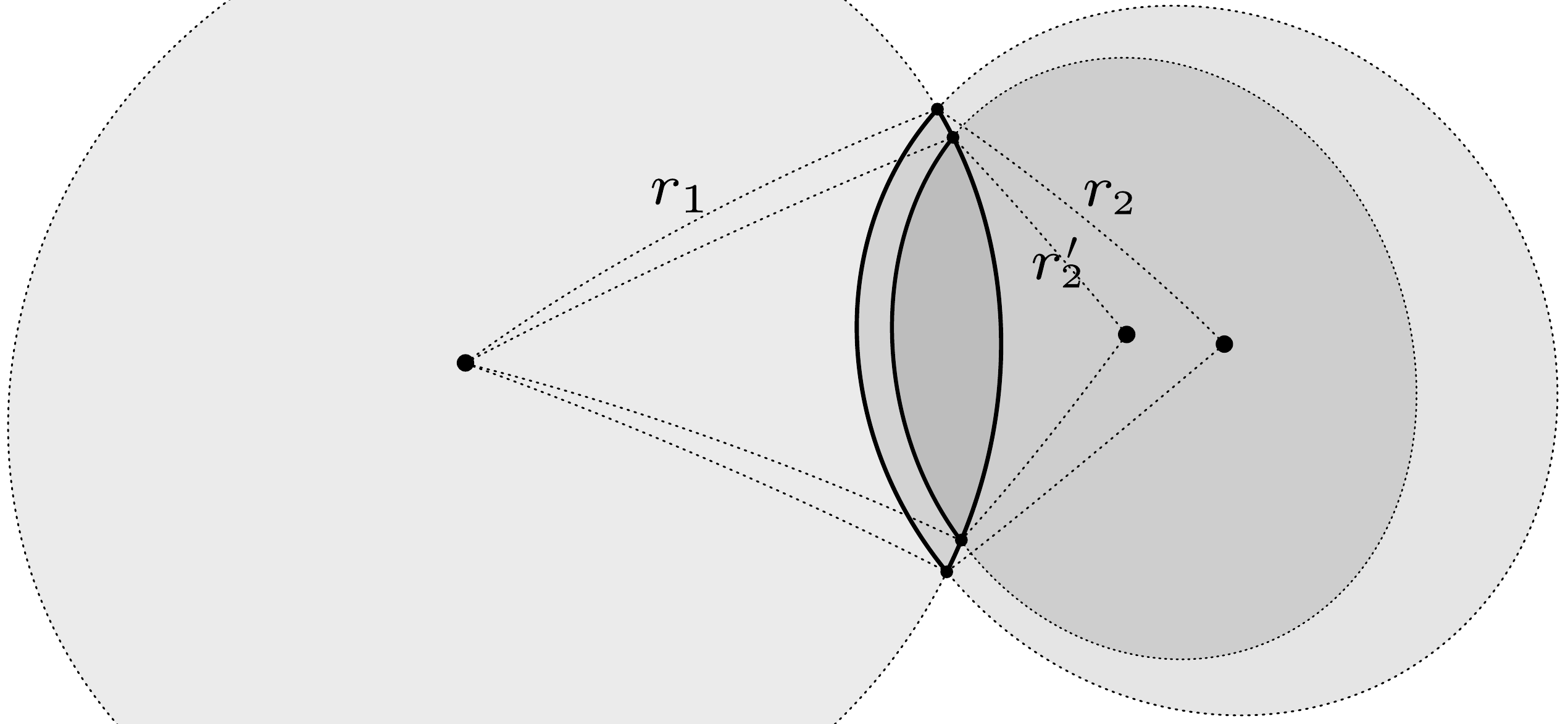}
	\caption{Variation of Figure \ref{figure_triangle} when $r_1$ is fixed and $r_2$ decreases.}
	\label{figure_shrink}
\end{figure}

\begin{lemma}\label{lemma_positive}
Consider $(K_1,K_2)$ as coordinates of $\mathcal{B}_\theta$. Then the matrix-valued function
$$
\begin{pmatrix}
\frac{\pa T_1}{\pa K_1}&\frac{\pa T_1}{\pa K_2}\\[3pt]
\frac{\pa T_2}{\pa K_1}&\frac{\pa T_2}{\pa K_2}
\end{pmatrix}
$$
(which is symmetric by Lemma \ref{lemma_closed}) is positively definition.
\end{lemma}
\begin{proof}
Using the symmetry $\frac{\pa T_1}{\pa K_2}=\frac{\pa T_2}{\pa K_1}$, for any $(x_1,x_2)\in\R^2$ we get
\begin{align*}
(x_1,x_2)
\begin{pmatrix}
	\frac{\pa T_1}{\pa K_1}&\frac{\pa T_1}{\pa K_2}\\[3pt]
	\frac{\pa T_2}{\pa K_1}&\frac{\pa T_2}{\pa K_2}
\end{pmatrix}
\binom{x_1}{x_2}&=-\frac{\pa T_1}{\pa K_2}(x_1-x_2)^2+\frac{\pa(T_1+T_2)}{\pa K_1}x_1^2+\frac{\pa(T_1+T_2)}{\pa K_2}x_2^2\\
&=-\frac{\pa T_1}{\pa K_2}(x_1-x_2)^2-\frac{\pa A}{\pa K_1}x_1^2-\frac{\pa A}{\pa K_2}x_2^2,
\end{align*}
where $A=2\theta-T_1-T_2$ is the area. When $K_1$ is fixed and $K_2$ increases (or equivalently, $r_1$ is fixed and $r_2$ decreases),
the bigon shrinks as in Figure \ref{figure_shrink}, so we have 
$$
\frac{\pa T_1}{\pa K_2}<0,\quad \frac{\pa A}{\pa K_2}<0.
$$
It also holds that $\frac{\pa A}{\pa K_1}<0$ for the same reason. The required statement follows.
\end{proof}

By Lemmas \ref{lemma_closed} and \ref{lemma_positive}, any primitive function $\mathcal{B}_\theta\to\R$ of the $1$-form $\omega_\theta$ is strictly convex in the coordinates $(K_1,K_2)$. We henceforth fix such a primitive, say
$$
\Omega_\theta(K_1,K_2):=\int_0^{(K_1,K_2)}\omega_\theta,
$$
and will use it to construct a function on $\mathcal{M}$. Note that $\Omega_\theta(K_1,K_2)$ is symmetric in $K_1$ and $K_2$ because $\omega_\theta$ is preserved by the involution $\mathcal{B}_\theta\to\mathcal{B}_\theta$ switching the labels 1 and 2.

\begin{remark}\label{remark_comparison}
On the space $\mathcal{B}_\theta^\mathbb{H}$ of hyperbolic bigons of angle $\theta$, the quantities listed at the beginning of this section still make sense, except that they are now related by $k_i=\coth r_i$ and $\ell_i=\alpha_i\sinh r_i$. The work of Bobenko and Springborn \cite{Bobenko-Springborn} explained in \S \ref{sec_1} is based on the primitive function of a $1$-form on $\mathcal{B}_\theta^\mathbb{H}$ with exactly the same expression as ours, which can also be written as
$$
\ell_1\dif k_1+\ell_2\dif k_2=-\frac{\alpha_1\dif r_1}{\sinh r_1}-\frac{\alpha_2\,\dif r_2}{\sinh r_2}=\alpha_1\dif u_1+\alpha_2\dif u_2,\ \ \text{where } u_i:=\int_{r_i}^{+\infty}\hspace{-0.2cm}\frac{\dif r_i}{\sinh r_i}=-\log \tanh\tfrac{r_i}{2}.
$$
This function is analogous to the one of Colin de Verdi\`ere \cite{ColindeVerdiere} on hyperbolic geodesic triangles and is related to $3$-dimensional hyperbolic volumes (see e.g.\ \cite{Luo_I}). 
A similar proof as that of Lemma \ref{lemma_positive} shows that the function is strictly convex with respect to the above variables $(u_1,u_2)$ rather than our variables $(K_1,K_2)$. The current paper is essentially the spherical counterpart of \cite{Bobenko-Springborn} with $u_i$ and $\alpha_i$ replaced by $K_i$ and $T_i$, respectively. 
\end{remark}

We will need the following fact on the total curvature parameters $T_1$ and $T_2$:
\begin{lemma}\label{lemma_T1T2}
For any $T_1,T_2>0$ with $T_1+T_2<2\theta$, there exists, up to isometry, a unique bigon of given angle $\theta$ whose two sides have total geodesic curvatures $T_1$ and $T_2$ respectively.
\end{lemma}
\begin{proof}
An equivalent statement is that the map $(T_1,T_2):\mathcal{B}_\theta\to\R_+^2$ is a bijection from $\mathcal{B}_\theta$ to the triangle $\Delta:=\{(T_1,T_2)\in\R_+^2\mid  T_1+T_2<2\theta\}$. By the above discussions, this map is the gradient of the strictly convex function $\Omega_\theta:\mathcal{B}_\theta\cong\R^2\to\R$ and has image contained in $\Delta$. But it is a basic fact that the gradient of a strictly convex $C^1$-function on $\R^n$ is a homeomorphism from $\R^n$ to a convex domain in $\R^n$. Therefore, our map is bijective from $\mathcal{B}_\theta$ to some convex subdomain $\Delta'$ of $\Delta$. We only need to show that $\Delta'=\Delta$, or equivalently, that the vertices $(0,0)$, $(2\theta,0)$ and $(0,2\theta)$ of $\Delta$ are limit points of the map. 

To this end, consider the sequences $a_n:=(-n,-n)$ and $b_n:=(n,-n)$ ($n=1,2,\cdots$) in $\mathcal{B}_\theta\cong\R^2$ (recall that we are using the $(K_1,K_2)$-parametrization of $\mathcal{B}_\theta$). For the bigon $a_n$, since the geodesic curvature $k_1(a_n)=k_2(a_n)=e^{-n}$ tends to $0$ as $n\to\infty$, the total curvature $T_i(a_n)=\ell_i(a_n)k_i(a_n)\leq\pi k_i(a_n)$ tends to $0$ as well, hence $(0,0)$ is the limit point of the image of $a_n$. For $b_n$, we have $T_2(b_n)\to0$ for the same reason, whereas $r_1(b_n)\to0$ because $\log\cot r_1(b_n)=K_1(b_n)=n\to +\infty$. It follows that the area $A(b_n)$, which is less than the area of a disk of radius $r_1(b_n)$, tends to $0$. By Gauss-Bonnet, we have $T_2(b_n)=2\theta-T_1(b_n)-A(b_n)\to2\theta$. Therefore, $(2\theta,0)$ is the limit point of the image of $b_n$. By symmetry, $(0,2\theta)$ is a limit point as well. 
\end{proof}
\begin{remark}\label{remark_hyperbolicbigon}
For hyperbolic bigons, the same statement holds for the angle parameters $\alpha_1$ and $\alpha_2$ in place of $T_1$ and $T_2$. This fits into the analogy between our setting and the hyperbolic one explained in Remark \ref{remark_comparison}.
\end{remark}

\section{Proof of Theorem \ref{thm_main}}
The rest of the proof is standard. We construct a function on the space $\mathcal{M}\cong \R^\mathsf{F}$ from \S \ref{sec_2} by adding up functions of the form $\Omega_{\theta_e}(K_{f(e^+)},K_{f(e^-)})$, where $K=(K_f)_{f\in \mathsf{F}}\in\R^\mathsf{F}$ is the coordinates, $e\in \mathsf{E}$ is an edge, and $f(e^\pm)\in \mathsf{F}$ are the faces incident with $e$ from the two sides (see the proof of Prop.\  \ref{prop_reconstruction} for notation). In other words, we add up the values of $\Omega_\theta$ at all the bigons of a $(\sigma,\mathcal{D})\in\mathcal{M}$. The sum is a strictly convex function because $\Omega_\theta$ is strictly convex on $\mathcal{B}_\theta$. 
A linear modification of this function solves the problem:

\begin{proposition}\label{prop_final}
Let $S$ be a closed surface and $\mathsf{G}$ be a graph weighted by $(\theta_e)_{e\in\mathsf{E}}$ as in Theorem \ref{thm_main}. Given $(T_f)_{f\in \mathsf{F}}\in\R_+^\mathsf{F}$, the following statements hold for the strictly convex function 
$$
\Omega:\mathcal{M}\cong\R^\mathsf{F}\to \R,\quad \Omega(K):=\sum_{e\in \mathsf{E}} \Omega_{\theta_e}(K_{f(e^+)},K_{f(e^-)})-\sum_{f\in \mathsf{F}}T_fK_f.
$$
\begin{enumerate}[label=(\roman*)]
\item\label{item_final1} $K$ is a critical point of $\Omega$ if and only if the corresponding $(\sigma,\mathcal{D})\in\mathcal{M}$ satisfies $T(\pa D_f)=T_f$ for all $f\in \mathsf{F}$.
\item\label{item_final2} $\Omega$ is proper if and only if $(T_f)$ satisfies 
\begin{equation}\label{eqn_Tconstraint}
	\sum_{f\in \mathsf{F}'}T_f<\sum_{
		\scalebox{0.65}{\parbox{2.9cm}{$e$: edge incident\\ with any face in $\mathsf{F}'$}}
	}\hspace{-0.5cm}2\theta_e\quad \text{for any subset of faces $\mathsf{F}'\subset \mathsf{F}$}.
\end{equation}
\end{enumerate}
\end{proposition}
\begin{proof}
\ref{item_final1} We use the notation from the proof of Prop.\ \ref{prop_reconstruction} and consider the triangulation therein, whose set of triangles $\hat{\mathsf{E}}$ identifies with the oriented edges of $\mathsf{G}$. Given $K\in\mathcal{M}\cong\R^\mathsf{F}$, each triangle $t\in\hat{\mathsf{E}}$ naturally corresponds to a spherical disk sector: in fact, the quadrilateral in Figure \ref{figure_reconstruction} is the union of the two sectors corresponding to the triangles $e^+$ and $e^-$,  respectively. Let $\hat{T}_t(K)$ denote the total geodesic curvature of the boundary arc of that sector.

By definition of the $1$-form $\omega_\theta=\dif \Omega_\theta$ in \S \ref{sec_3}, we have
\begin{equation}\label{eqn_finalproof1}
\dif \sum_{e\in \mathsf{E}} \Omega_{\theta_e}(K_{f(e^+)},K_{f(e^-)})=\sum_{e\in \mathsf{E}}\Big(T_1(K_{f(e^+)},K_{f(e^-)})\dif K_{f(e^+)}+T_2(K_{f(e^+)},K_{f(e^-)})\dif K_{f(e^-)}\Big).
\end{equation}
But with the above notation, we may write $T_1(K_{f(e^+)},K_{f(e^-)})=\hat{T}_{e^+}(K)$ and $T_2(K_{f(e^+)},K_{f(e^-)})=\hat{T}_{e^-}(K)$. Thus, we can rewrite \eqref{eqn_finalproof1} as a sum over all triangles $t\in\hat{\mathsf{E}}$ and get
$$
\dif \sum_{e\in \mathsf{E}} \Omega_{\theta_e}(K_{f(e^+)},K_{f(e^-)})=\sum_{t\in \hat{\mathsf{E}}}\hat{T}_t(K)\dif K_{f(t)}=\sum_{f\in \mathsf{F}}T(\pa D_f)\dif K_f,
$$
where the last equality is obtained by first summing over those $t\in\hat{\mathsf{E}}$ contained in a face $f$. The required statement follows from the last expression.

\ref{item_final2} This is base on the fact, proved in \cite[Prop.\ 4]{Bobenko-Springborn} using linear programming (see also \cite{Guo_note,Rivin}), that the condition ``$(T_f)$ satisfies \eqref{eqn_Tconstraint}'' is equivalent to the existence of a  \emph{coherent angle system}\footnote{
In \cite{Bobenko-Springborn, ColindeVerdiere}, such a system represents angles. Here we deal with a system defined by the same conditions, hence keep the name, although it now represents total geodesic curvatures rather angles (cf.\ Remark \ref{remark_hyperbolicbigon}).
}
for $(T_f)$, which is by definition some $(\hat{T}_t)_{t\in\hat{\mathsf{E}}}\in\R_+^{\hat{\mathsf{E}}}$ such that for any edge $e$ and any face $f$ of $\mathsf{G}$, we have
$$
\hat{T}_{e^+}+\hat{T}_{e^-}<2\theta_e,\quad \sum_{\scalebox{0.65}{$t$: triangle in $f$}}\hspace{-0.2cm}\hat{T}_t=T_f.
$$ 

By this fact, it suffices to show that $\Omega$ is proper if and only if there is a coherent angle system $(\hat{T}_t)$ for $(T_f)$. 

If $\Omega$ is proper, then $\Omega$ has a critical point $K\in\mathcal{M}\cong\R^\mathsf{F}$. By Statement \ref{item_final1}, the above defined $(\hat{T}_t(K))_{t\in\hat{\mathsf{E}}}$ is a coherent angle system  for $(T_f)$.  This proves the ``only if'' part.

Conversely, given a coherent angle system $(\hat{T}_t)$  for $(T_f)$, we may rewrite the linear part of $\Omega(K)$ as 
$$
\sum_{f\in\mathsf{F}}T_fK_f=\sum_{e\in\mathsf{E}}\left(\hat{T}_{e^+}K_{f(e^+)}+\hat{T}_{e^-}K_{f(e^-)}\right),
$$
and thus write
\begin{equation}\label{eqn_Omega}
\Omega(K)=\sum_{e\in\mathsf{E}}\left(\Omega_{\theta_e}(K_{f(e^+)},K_{f(e^-)})-\hat{T}_{e^+}K_{f(e^+)}-\hat{T}_{e^-}K_{f(e^-)}\right).
\end{equation}
Each summand in \eqref{eqn_Omega} can be understood as a function on the bigon space $\mathcal{B}_\theta\cong \R^2$ of the form
$$
\Omega_{\theta,\hat{T}_1,\hat{T}_2}(K_1,K_2):=\Omega_\theta(K_1,K_2)-\hat{T}_1K_1-\hat{T}_2K_2
$$ 
for constants $\hat{T}_1,\hat{T}_2>0$, $\hat{T}_1+\hat{T}_2<2\theta$. The differential of this function is
$$
\dif\Omega_{\theta,\hat{T}_1,\hat{T}_2}=\omega_\theta-\hat{T}_1\dif K_1-\hat{T}_2\dif K_2
=(T_1-\hat{T}_1)\dif K_1+(T_2-\hat{T}_2)\dif K_2.
$$
Therefore, $\Omega_{\theta,\hat{T}_1,\hat{T}_2}$ has a critical point at the bigon whose two sides have total geodesic curvatures $\hat{T}_1$ and $\hat{T}_2$, which exists by Lemma \ref{lemma_T1T2}. Since a strictly convex function defined on the whole $\R^n$ is proper if and only if it has a critical point, we infer that $\Omega_{\theta,\hat{T}_1,\hat{T}_2}$ is proper on $\mathcal{B}_\theta$, and hence the whole sum \eqref{eqn_Omega} is proper on $\mathcal{M}$. This proves the ``if'' part.
\end{proof}

To deduce Theorem \ref{thm_main}, note that by Proposition \ref{prop_final} \ref{item_final1}, we may understand the theorem as saying that $(T_f)\in\R_+^\mathsf{F}$ satisfies \eqref{eqn_Tconstraint} if and only if $\Omega$ has a critical point, and that the critical point is unique if it exists. The last uniqueness follows from the strict convexity, while the ``if and only if'' statement follows from Proposition \ref{prop_final} \ref{item_final2} and the fact about strictly convex functions quoted at the end of the above proof.

\bibliographystyle{siam} \bibliography{pattern}
\end{document}